\documentclass[12pt]{article}
\usepackage{amsmath, graphicx}
\usepackage{amsthm}
\usepackage[english]{babel}
\usepackage{amsfonts, amssymb}
\newtheorem{theorem}{Theorem}
\newtheorem{corollary}[theorem]{Corollary}
\newtheorem{lemma}[theorem]{Lemma}
\newtheorem{proposition}[theorem]{Proposition}
\theoremstyle{remark}
\newtheorem{definition}[theorem]{Definition}

\newtheorem{example}{Example}

\begin{document}

%%%%%%%%%%%%%%%%%%%%%%%%%%%%%%%%%%%%%%%%%%%%%%%%%%%
%%%%                                           %%%%
%%%%   Generic identifiability and 2nd order   %%%%
%%%%  sufficiency in tame convex optimization  %%%%
%%%%                                           %%%%
%%%%   J. Bolte, A. Daniilidis, A. Lewis       %%%%
%%%%                                           %%%%
%%%%              January 18, 2009             %%%%
%%%%             [Version: arisd]              %%%%
%%%%                                           %%%%
%%%%%%%%%%%%%%%%%%%%%%%%%%%%%%%%%%%%%%%%%%%%%%%%%%%

\title{Generic identifiability and second-order sufficiency in tame convex optimization}
\author{J. Bolte \thanks{This author acknowledges the support of the grant ANR-05-BLAN-0248-01
(France).}, A. Daniilidis \thanks{This author's research
was partially supported by the grants ANR-05-BLAN-0248-01 (France)
and MTM2008-06695-C03-03 (Spain)}, and A.~S.~Lewis
\thanks{This author's research was partially supported by the ANR grant
ANR-05-BLAN-0248-01 (France) and by the National Science Foundation
Grant DMS-0806057 (USA)}}

\date{Date:  \today~~ }

\maketitle

\begin{abstract}
We consider linear optimization over a fixed compact convex feasible region
that is semi-algebraic (or, more generally, ``tame'').  Generically, we prove
that the optimal solution is unique and lies on a unique manifold, around which
the feasible region is ``partly smooth'',  ensuring finite identification of the
manifold by many optimization algorithms.  Furthermore, second-order optimality
conditions hold, guaranteeing smooth behavior of the optimal solution under small
perturbations to the objective.
\end{abstract}

%\bigskip


\bigskip

\textit{Secondary} 32B20, 90C05, 49K40

\section{Introduction}
``Identification'' in constrained optimization signifies an
important idea both in theory and for algorithms. Sensitivity
analysis, the theory of how optimal solutions behave under data
perturbations, depends on identifying active constraints and
verifying associated optimality conditions. A variety of practical
algorithms for inequality-constrained problems aim to identify the
active constraints:  once the identification is successful, we have
essentially converted to the easier, equality-constrained case.
\smallskip

An early survey of identification techniques, for optimization over polyhedra and generalizations, appears in \cite{Bur88}. For general convex
feasible regions, a more abstract approach appeals, in part for its
theoretical elegance, and in part because for constraints more
complex than simple inequalities, such as the semidefinite
inequalities common in modern optimization, simply deciding whether
a constraint is active or not fails to capture crucial finer
details.  Such an abstract approach, based on the idea of an
``identifiable surface'', appeared in \cite{Wri93}.  As shown in
\cite{Lewis2002}, this idea has an equivalent but more geometric
description: the surface turns out to be a manifold contained in the
feasible region satisfying a property called ``partial smoothness''.
As well as its geometric transparency, the notion of partial
smoothness has the merit of extending naturally to the nonconvex
case.  In this work, however, we confine ourselves to convex
feasible regions.
\smallskip

Our goal here is to show that partial smoothness is a common
phenomenon. Certainly the property can fail, either because the
feasible region is somehow pathological, or because of the failure
of the typical regularity conditions needed for standard sensitivity
analysis. A good illustration is the convex optimization problem
over $\mathbb{R}^3$,
\begin{equation} \label{bad example}
\inf\{ w : w \ge (|u| + |v|)^2 \}.
\end{equation}
As we perturb the linear objective function slightly, the corresponding
optimal solutions describe not one but two distinct manifolds.
\smallskip

Nonetheless, very generally, partial smoothness is indeed typical
for linear optimization over a fixed compact convex feasible region
$F \subset\mathbb{R}^n$. Specifically, we prove that if $F$ is
\emph{semi-algebraic}---a finite union of sets defined by finitely
many polynomial inequalities---or, more generally, ``tame'', then,
except for objectives lying in some exceptional set of dimension
strictly smaller than $n$, the corresponding optimal solution is
unique, and $F$ is partly smooth around a corresponding unique
manifold.  Furthermore, a second-order sufficient optimality
condition holds.  A variety of algorithms will therefore identify
the manifold and converge well, and standard sensitivity analysis
applies.
\smallskip

Various authors have shown that, for some suitably structured convex
optimization problem, the set of instances for which the optimal
solution has some beneficial property (such as
``well-posedness''~\cite{DonZol93}) is generic. An interesting
recent example is~\cite{IoffeLuc2005}. By contrast, we assume
nothing about our feasible region, beyond its semi-algebraic or tame
nature. Remarkably, nonetheless, generic problems are very well
behaved.

\section{Preliminaries and notation}
Throughout the manuscript we deal with a finite-dimensional
Euclidean space $\mathbb{R}^{n}$ equipped with the usual scalar
product $\langle\cdot,\cdot\rangle$ and the corresponding Euclidean
norm $\Vert\cdot\Vert$. We denote by $B(x,r)$ the closed ball with
center $x\in\mathbb{R}^{n}$ and radius $r>0.$ We simply denote by
$\mathbb{B}$ the closed unit ball $B(0,1)$ and by $S^{n-1}$ its
boundary, that is, the unit sphere of $\mathbb{R}^{n}$. Given any $E
\subset \mathbb{R}^n$, we denote by $\mathrm{ri\,}E$ its relative
interior and by $\overline{E}$ its closure.

\subsubsection*{Preliminaries on variational analysis}
We refer to \cite{Phelps93} and \cite{Rock98} for basic facts about convex and variational
analysis that we use.
\smallskip

\noindent Let $X,Y$ be metric spaces and $T:X\rightrightarrows Y$ be
a set-valued mapping. We say that $T$ is \emph{outer semicontinuous}
at a point $\bar x \in X$ if, for any sequence of points $x_r \in X$
converging to $\bar x$ and any sequence of points $y_r \in T(x_r)$
converging to $\bar y$, we must have $\bar y \in T(\bar x)$. On the
other hand, we say that $T$ is \emph{inner semicontinuous} at $\bar
x$ if, for any sequence of points $x_r \in X$ converging to $\bar x$
and any point $\bar y \in Y$, there exists a sequence $y_r \in Y$
converging to $\bar y$ such that $y_r \in T(x_r)$ for all large $r$.
If both properties hold, we call $T$ \emph{continuous} at $\bar x$.
\smallskip

\noindent Consider a nonempty closed convex set $F \subset
\mathbb{R}^n$. The \emph{normal cone} $N_{F}(x)$ at a point $x \in
F$ is defined as follows:
\begin{equation}
N_{F}(x) ~=~ \Big\{c\in\mathbb{R}^{n} : \langle
c,x^{\prime}-x\rangle\leq0,~\forall x^{\prime}\in F \Big\}.
\label{normal cone}
\end{equation}
It is standard and easy to check that the mapping $x \mapsto
N_{F}(x)$ is outer semicontinuous on $F$. In a slightly different
context, for a point $x$ in a smooth submanifold ${\mathcal M}$ of
$\mathbb{R}^n$, we denote by $N_{\mathcal M}(x)$ the normal space in
the usual sense of elementary differential geometry, that is, the
orthogonal complement in $\mathbb{R}^n$ of the tangent space
$T_{\mathcal{M}}(x)$ of $\mathcal{M}$ at $x$.

\subsubsection*{Preliminaries on partial smoothness}
We recall from \cite{Lewis2002} the definition of partial smoothness, specialized to the convex case.

\begin{definition}
\label{Definition_partial-smoothness} A closed convex set $F \subset
\mathbb{R}^n$ is called \emph{partly smooth} at a point $\bar{x} \in
F$ relative to a set $\mathcal{M} \subset F$ if the following
properties hold:
\begin{enumerate}
\item[(i)]
$\mathcal{M}$ is a $C^{2}$ submanifold of $\mathbb{R}^{n}$ (called the \emph{active} manifold) containing $\bar{x}$.
\item[(ii)]
The set-valued mapping $x \mapsto N_F(x)$, restricted to the domain $\mathcal{M}$, is continuous at $\bar{x}$.
\item[(iii)]
$
N_{\mathcal{M}}(\bar{x}) = N_{F}(\bar{x})-N_{F}(\bar{x}).
$
\end{enumerate}
\end{definition}

\noindent While not obvious from the above definition, the active
manifold for a partly smooth convex set is locally unique around the
point of interest:  see \cite[Cor.~4.2]{Lewis2004}.
\smallskip

Geometrically, condition (iii) guarantees ``sharpness'' around a
kind of ``ridge'' in the set $F$ defined by the active manifold, as
illustrated in the following simple example.
\begin{example}
In $\mathbb{R}^3$, define
\begin{eqnarray*}
F & = & \{(u,v,w) : w\geq u^{2}+|v|\},\\
\mathcal{M} & = & \{(t,0,t^2) : t\in(-1,1)\}.
\end{eqnarray*}
Then the set $F$ is partly smooth at the point $\bar{x}=(0,0,0)$ relative to
the one-dimensional manifold $\mathcal{M}$.
\end{example}
\smallskip

The following example illustrates the importance of normal cone continuity.

\begin{example}
[failure of normal cone continuity]
In $\mathbb{R}^3$, consider the set and manifold
\begin{eqnarray*}
F & = & \Big\{ (u,v,w) : v\geq0,~w\geq0,~v+w\geq u^{2} \Big\} \\
\mathcal{M} & = & \{(t,t^{2},0) : t \in (-1,1)\}.
\end{eqnarray*}
Then $F$ is convex and conditions (i) and (iii) of
Definition~\ref{Definition_partial-smoothness} are satisfied at the
point $\bar x=(0,0,0)$. But condition (ii) fails, since the normal
cone mapping is discontinuous there, relative to $\mathcal{M}$.
\end{example}

For the purposes of sensitivity analysis, partial smoothness is
most useful when combined with a second-order sufficiency condition,
captured by the following definition.

\begin{definition}
\label{Definition_strong-critical} Consider a vector $\bar{c} \in
\mathbb{R}^n$ and a closed convex set $F \subset \mathbb{R}^n$ that
is partly smooth at a point $\bar{x} \in \mbox{argmax}_F
\langle\bar{c},\cdot\rangle$ relative to a manifold $\mathcal{M}$.
We say that $\bar{x}$ is \emph{strongly critical} if the following
properties hold:
\begin{enumerate}
\item[(i)]
$\bar{c}\in\mathrm{ri\,}N_{F}(\bar{x})$.
\item[(ii)]
There exists $\delta>0$ such that
\[
\langle\bar{c},\bar{x}\rangle\geq\langle\bar{c},x\rangle+\delta||x-\bar
{x}||^{2}\text{,\quad for all }x\in\mathcal{M}\;\text{near }\bar{x}.
\]
\end{enumerate}
\end{definition}
Condition (i) can be interpreted as a kind of ``strict
complementarity'' condition, while condition (ii) concerns quadratic
decay. Notice that the above definition yields uniqueness of the
maximizer $\bar{x}$ of $\bar{c}$, as well as good sensitivity
properties, as the following result shows: see \cite{Lewis2002}.

\begin{theorem}[second-order sufficiency]
Consider a closed convex set $F \subset \mathbb{R}^n$ and assume
that $F$ is partly smooth at some point $\bar{x}\in F$ and that
$\bar{x}$ is strongly critical point  for the problem $\max_F
\langle\bar{c},\cdot\rangle$, relative to a manifold $\mathcal{M}$.
Then for all vectors $c \in \mathbb{R}^n$ sufficiently near
$\bar{c}$, the perturbed problem $\max_F \langle c,\cdot\rangle$ has
a unique optimal solution $x_c \in \mathcal{M}$. The map $c \mapsto
x_c$ is $C^1$ around $\bar{c}$.
\end{theorem}

\subsubsection*{Preliminaries on tame geometry}
Let us first recall the definitions of an ``o-minimal structure''
(see for instance \cite{Coste99}, \cite{Dries-Miller96} or
\cite{Kurdyka98} and references therein). \smallskip

\begin{definition}
\label{Definition_o-minimal} An \emph{o-minimal structure} on
$(\mathbb{R},+,.)$ is a sequence of Boolean algebras
$\mathcal{O}=\{\mathcal{O}_{n}\}$, where each algebra
$\mathcal{O}_n$ consists of subsets of $\mathbb{R}^{n}$, called
\emph{definable} (in $\mathcal{O}$), and such that for every
dimension $n\in\mathbb{N}$ the following properties hold.
\begin{enumerate}
\item[(i)]
For any set $A$ belonging to $\mathcal{O}_{n}$, both $A\times\mathbb{R}$ and
$\mathbb{R}\times A$ belong to $\mathcal{O}_{n+1}$.
\item[(ii)]
If $\Pi:\mathbb{R}^{n+1}\rightarrow\mathbb{R}^{n}$ denotes the
canonical projection, then for any set $A$ belonging to
$\mathcal{O}_{n+1}$, the set $\Pi(A)$ belongs to $\mathcal{O}_{n}$.
\item[(iii)]
$\mathcal{O}_{n}$ contains every set of the form $\{x\in\mathbb{R}^{n} : p(x)=0\}$,
for polynomials $p:\mathbb{R}^{n}\rightarrow\mathbb{R}$.
\item[(iv)]
The elements of $\mathcal{O}_{1}$ are exactly the finite unions of
intervals and points.
\end{enumerate}
When $\cal{O}$ is a given o-minimal structure, a function
$f:\mathbb{R}^{n}\rightarrow\mathbb{R}^{m}$ (or a set-valued mapping
$F:\mathbb{R}^{n}\rightrightarrows\mathbb{R}^{m}$) is called \emph{definable}
(in $\mathcal{O}$) if its graph is definable as a subset of $\mathbb{R}^{n}\times\mathbb{R}^{m}$.
\end{definition}

If a subset $A$ of $\mathbb{R}^{n}$ has the property that its
intersection with every ball is definable in some o-minimal
structure, then it is sometimes called \emph{tame}. In this work we
are concerned primarily with bounded sets:  in that context, we use
the terms ``tame'' and ``definable'' interchangeably.
\smallskip

Semi-algebraic sets constitute an o-minimal structure, as a
consequence of the Tarski-Seidenberg principle, but richer
structures also exist. In particular, the Gabrielov theorem implies
that ``subanalytic'' sets are tame. These two structures in
particular provide rich practical tools, because checking
semi-algebraicity or subanalyticity of sets in concrete problems of
variational analysis is often easy. We refer to \cite{BDL04},
\cite{BDL2008}, and \cite{Ioffe2008} for more details. \smallskip

Definable sets and functions enjoy many structural properties. In
particular, every definable set can be written as a finite disjoint
union of manifolds (or ``strata'') that fit together in a regular
``stratification'': see \cite[\S 4.2]{Dries-Miller96}. In
particular, the \emph{dimension} of the set is the maximum of the
dimensions of the strata, a number independent of the
stratification:  see \cite[Definition~9.14]{Coste99} for more
details. We call a definable subset of a definable set
\emph{generic} if its complement has strictly smaller dimension.
\smallskip

In this paper we make fundamental use of a stratification result. We
present a particular case---adapted to our needs---of a more general
result: see \cite[p.\ 502, \S 1.19 (2)]{Dries-Miller96} or
\cite{TaLeLoi96} for the statement in its full generality.  The
result describes a decomposition of the domain of a definable
function into subdomains on which the function has ``constant
rank'':  a smooth function has \emph{constant rank} if its
derivative has constant rank throughout its domain.  Such functions
have a simple canonical form:  they are locally equivalent to
projections, as described by the following result from basic
differential geometry (see \cite[Thm 7.8]{Lee2003}.

\begin{proposition}
[Constant Rank Theorem]\label{Proposition_constant}
Let $M_1$ and $M_2$ be two differentiable manifolds, of dimensions
$m_1$ and $m_2$ respectively, and let $g:M_1\rightarrow M_2$ be a
differentiable mapping of constant rank $r$. Then for every point
$x\in M_1$, there exist neighborhoods $O_i$ of zero in $\mathbb{R}^{m_i}$
and local diffeomorphisms $\psi_i :O_i\rightarrow M_i$ (for $i=1,2$)
with $\psi_1(0)=x$ and $\psi_2(0)=g(x)$, such
that mapping $\psi_2^{-1} \circ g \circ \psi_1$ is just the projection
$\pi : O_1 \rightarrow O_2$ defined by
\begin{equation}
\pi(y_1,y_2,\ldots,y_{m_{1}})
~=~
(y_1,y_2,\ldots,y_r,0,\ldots,0)
\in\mathbb{R}^{m_2},~~ (y \in O_1).
\label{projection}
\end{equation}
\end{proposition}

The stratification result we use follows.

\begin{proposition}
[Constant rank stratification]
\label{Proposition_rank}
Let $f:M\rightarrow \mathbb{R}^{n}$ be a definable function, where $M$
is a submanifold of $\mathbb{R}^{n}$. Then there exists a $C^2$-stratification
$\mathcal{S}=\{S_{i}\}_{i}$ of $M$ and a $C^2$-stratification
$\mathcal{T}$ of $\mathbb{R}^{n}$ such that the restriction $f_{i}$
of $f$ onto each stratum $S_{i}\in\mathcal{S}$ is a
$C^2$-function, $f_{i}(S_{i})\in\mathcal{T}$ and $f_{i}$ is of
constant rank in $S_{i}.$
\end{proposition}

The above statement yields that each restriction
$f_{i}:S_{i}\rightarrow$ $f_{i}(S_{i})$ is surjective, $C^2$, and of
constant rank $r_{i}$. Thus $r_{i}$ is also equal to the dimension
of the manifold $f_{i}(S_{i})$:
\[
r_{i}=\mathrm{rank\,}f_{i}=\dim\operatorname{Im}(df_{i}(x))=\dim
(T_{f_{i}(S_{i})})(f_{i}(x)),~~ \text{for all}~ x\in S_{i}.
\]

\section{Introductory results}
We always consider a fixed nonempty compact convex set $F \subset \mathbb{R}^n$,
and study the set of optimal solutions of the problem
\[
\sup_F \langle c, \cdot \rangle
\]
for vectors $c \in \mathbb{R}^n$. The optimal value of this problem,
as a function of $c$, is called the \emph{support function}, denoted
$\sigma_F$. We denote by $\mbox{argmax}_F \langle c ,\cdot\rangle$
the set of optimal solutions. By scaling, we may as well assume $c$
lies in the unit sphere $S^{n-1}$. We aim to show good behavior for
objective vectors $c$ lying in some large subset of the sphere.
Classically, ``large'' might mean, for example, ``full-measure'', or
perhaps ``generic'':  a \emph{generic} subset of a topological space
is one containing a countable intersection of dense open sets.
Clearly, these distinctions vanish for definable sets.
\smallskip

We begin our development with an easy and standard argument.

\begin{proposition}[Generic uniqueness]
Consider a nonempty compact convex set $F \subset \mathbb{R}^n$. For
all vectors $c$ lying in a generic and full-measure subset of the
sphere $S^{n-1}$, the linear functional $\langle c,\cdot \rangle$
has a unique maximizer over $F$.
\end{proposition}

\noindent\textbf{Proof.} We use various standard techniques from
convex analysis \cite{Phelps93,Rock98}. Note first
\[
\mbox{argmax}_F \langle c,\cdot\rangle = \partial\sigma_{F}(c),
\]
where $\partial$ denotes the convex subdifferential. This set is
therefore a singleton if and only if the support function $\sigma_F$
is differentiable at $c$.  Being a finite, positively homogeneous,
convex function, the set of points of differentiability is both
generic and full-measure in $\mathbb{R}^n$, and is closed under
strictly positive scalar multiplication. The result now follows.
$\hfill\Box$
\bigskip

As a next step towards our main result, we prove stronger properties
for at least a dense set of objectives. Density will suffice for our
purposes once we move to a tame setting.

\begin{proposition}
[Almost all linear functionals have strong maximizers]
\label{Proposition_Lemma 5} Let $F$ be a nonempty compact convex
subset of $\mathbb{R}^{n}$. Then corresponding to any vector $c$
lying in some subset of $S^{n-1}$ of full measure, there exist a
vector $x_{c}\in F$ and a constant $\delta_{c}>0$ such that
\begin{equation}
\langle c,x_{c}\rangle\geq\langle
c,x\rangle+\delta_{c}||x-x_{c}||^{2} ,\quad\text{for all }x\in
F,\label{strong}
\end{equation}
that is, $x_{c}$ is a strong (unique) maximizer of the linear functional
$\langle c,\cdot\rangle$ over $F$.
\end{proposition}

\noindent\textbf{Proof.} Let us denote  by  $\sigma_{F}$ the support
function of $F$ and $i_{F}$ the corresponding indicator
function
\[
i_{F}(c):=\left\{
\begin{array}{l}0\mbox{ if }x\in F\\
+\infty \mbox{ otherwise.}
\end{array}\right.\]
Notice that $\sigma_{F}$, a finite convex function, coincides with
the Fenchel conjugate of $i_{F}$ and $\partial\sigma_{F}(c)$
is the set of maximizers of $c$ on $F$. Applying Alexandrov's
Theorem (\cite[Theorem 13.51, p. 626]{Rock98}), we deduce that there
exists a full measure subset $A$ of $\mathbb{R}^{n}$ on which
$\sigma_{F}$ has a quadratic expansion. (Thus in particular
$\sigma_{F}$ is differentiable there and $\nabla\sigma
_{F}(c)=x_{c}$, where $x_{c}$ denotes the unique maximizer of $c$
at $F$.) In view of \cite[Definition 13.1(c), p. 580]{Rock98}, we
have, for any fixed $\bar{c}\in A$, there exists a positive semidefinite matrix $S$ such that for all $c\in\mathbb{R}^n$,
\[
\sigma_{F}(c)=\sigma_{F}(\bar{c})+\langle\nabla\sigma_{F}(\bar{c}),c-\bar
{c}\rangle+\frac{1}{2}\langle
S(c-\bar{c}),c-\bar{c}\rangle+o(||c-\bar {c}||^{2})\,.
\]
Hence there exists $\varepsilon>0$, $\rho>0$ such that for all $c\in
B(\bar{c},\varepsilon)$ we have
\[
\sigma_{F}(c)\leq\sigma_{F}(\bar{c})+\langle x_{\bar{c}},c-\bar{c}\rangle
+\frac{\rho}{2}||c-\bar{c}||^{2}\,.
\]
Further, we can clearly assume
\begin{equation}
\varepsilon^{-1}\,\mathrm{diam}\,(F)\,<\,\rho.\label{ar1}
\end{equation}
Now consider $x\in F$. Recalling $\sigma_{F}(\bar{c})=\langle\bar
{x}_{c},\bar{c}\rangle$ we deduce successively
\begin{align*}
0=i_{F}(x)=\sigma_{F}^{\ast}(x)  & =\sup_{c\in\mathbb{R}^{n}}\left\{
\,\langle x,c\rangle-\sigma_{F}(c)\,\right\}  \\
& \geq\sup_{c\in B(\bar{c},\varepsilon)}\left\{  \,\langle x,c\rangle
-\sigma_{F}(c)\,\right\}  \\
& \geq\sup_{c\in B(\bar{c},\varepsilon)}\left\{  \langle x,c\rangle-\sigma
_{F}(\bar{c})-\langle x_{\bar{c}},c-\bar{c}\rangle-\frac{\rho}{2}||c-\bar
{c}||^{2}\right\}  \\
& =\sup_{c\in B(\bar{c},\varepsilon)}\left\{  \langle x-x_{\bar{c}},
c\rangle-\frac{\rho}{2}||c-\bar{c}||^{2}\right\}  \\
& =\langle x-x_{\bar{c}},\bar{c}\rangle+\sup_{u\in B(0,\varepsilon)}\left\{
\langle x-x_{\bar{c}},u\rangle-\frac{\rho}{2}||u||^{2}\right\}  .
\end{align*}
In view of \eqref{ar1} it is easy to notice that the above supremum
is realized at $u=\rho^{-1}(x-x_{\bar{c}})\in B(0,\varepsilon).$
Replacing this value in the above inequality we deduce
\[
0\geq\langle
x-x_{\bar{c}},\bar{c}\rangle+\frac{1}{2\rho}||x-x_{\bar{c}}
||^{2},\qquad\text{for all }x\in F.
\]
which yields the asserted equation for $\delta_{c}=(2\rho)^{-1}$.
The restriction of the result to $S^{n-1}$ is straightforward. $\hfill\Box$

\bigskip

\begin{corollary}[Density of functionals with strong maximizer]\label{dense}
The set of vectors $c\in S^{n-1}$ such that there exist $x_{c}\in
F$, $c\in\mbox{\rm ri}\,N_{F}(x_{c})$ and a constant $\delta_{c}>0$
such that
\[
\langle c,x_{c}\rangle\geq\langle
c,x\rangle+\delta_{c}||x-x_{c}||^{2} ,\quad\text{for all }x\in F,
\]
is a dense subset of the sphere $S^{n-1}$.
\end{corollary}

\noindent\textbf{Proof} Given $c\in A$, take $c^{\prime}$ in
$\mbox{\rm ri}\,N_{F}(x_{c})$ and choose $\eta>0$ such that $c^{\prime}-\eta
c\in\mbox{\rm ri}\,N_{F}(x_{c})$. From the definition of the normal cone we
deduce
\[
\langle c^{\prime},x_{c}-x\rangle=\langle\eta c,x_{c}-x\rangle+\langle
c^{\prime}-\eta c,x_{c}-x\rangle\geq\eta\delta_{c}||x-x_{c}||^{2}=\eta
\delta_{c}||x-x_{c^{\prime}}||^{2},
\]
for all $x$ in $F$. In other words: if $\mathcal{C}$ denotes the set
of linear functionals $c\in\mathbb{R}^{n}$ satisfying (\ref{strong})
with $c\in \mbox{\rm ri}\,N_{F}(x_{c})$, then
$\overline{\mathcal{C}}\supset A$ and thus
$\overline{\mathcal{C}}=\overline{A}=\mathbb{R}^{n}$. The density
result on $S^{n-1}$ follows easily. $\hfill\Box$

\section{Main result}
>From now on we shall assume that the nonempty compact convex set $F
\subset \mathbb{R}^n$ is also definable in some o-minimal structure
(see Definition~\ref{Definition_o-minimal}). We are ready to state
and prove the main result of this work. This result asserts that a
generic linear optimization problem over $F$ has a unique optimal
solution, that $F$ is partly smooth there, and strong criticality
holds. As we see in the proof below, the active manifold arises
naturally, by means of Proposition~\ref{Proposition_rank} (constant
rank stratification) applied to an appropriately defined function.
\smallskip

To obtain the semi-algebraic version of the result below, simply
replace the term ``definable'' by ``semi-algebraic''.

\begin{theorem}
[Main result]\label{Theorem_main} Let $F$ be a nonempty compact
convex subset of $\mathbb{R}^{n}$ that is definable in some
o-minimal structure. Then there exists a definable generic subset
$U$ of the unit sphere $S^{n-1}$ with the following property: for
each unit vector $c \in U$, there exists a unique vector $x_c \in F$
and a definable set $\mathcal{M}_{c} \subset F$ (unique in a
neighborhood of $x_c$)  satisfying:
\begin{enumerate}
\item[\mbox{\rm (i)}]
$\mbox{\rm argmax}_F \langle c, \cdot \rangle = \{x_{c}\}$;
\item[\mbox{\rm (ii)}]
$F$ is partly smooth at $x_c$ relative to $\mathcal{M} _{c}$;
\item[\mbox{\rm (iii)}]
$x_{c}$ is strongly critical.
\end{enumerate}
\end{theorem}

\noindent\textbf{Proof}
Let us consider the definable set-valued
mapping $\tilde{\Phi}:S^{n-1}\rightrightarrows F$ defined by
\begin{equation}
\tilde{\Phi}(c)= \mbox{argmax}_F \langle
c, \cdot \rangle,\label{a1}
\end{equation}
and let us note, by the definition of the normal cone,
\[
\tilde{\Phi}^{-1}(x)=N_{F}(x)\cap S^{n-1}.
\]
Let $D$ denote the dense subset of $S^{n-1}$ asserted in
Corollary~\ref{dense} (Density of functionals with
strong maximizer). Since $F$ is a definable set, we deduce easily
that $D$ is also definable (see \cite[Section~2.2]{BDL04}, for
example), and hence generic. In particular, the set
\begin{equation}
N_{\ast}=S^{n-1}\setminus D \label{a2}
\end{equation}
has dimension strictly less than $n-1$.
\smallskip

\noindent Let $\Phi:D\rightarrow F$ denote the restriction of the
mapping $\tilde{\Phi}$ to $D$. Observe that $\Phi$ is single-valued
and, by the definition of the set $D$, satisfies the strict
complementarity and quadratic decay conditions:
\begin{enumerate}
\item[(i)]
$c\in\mathrm{ri\,}N_{F}(\Phi(c))$;
\item[(ii)]
$\langle c,\Phi(c)\rangle\geq\langle
c,x\rangle+\delta||x-\Phi(c)||^{2} $,\quad for some $\delta>0$ and
for all $x\in F$.
\end{enumerate}
Applying Proposition~\ref{Proposition_rank} (Constant rank
stratification) to the definable function
\begin{equation}\label{aris-C1}
\left\{
\begin{array}{l}
 \Phi : D\rightarrow F \smallskip \\
 c\mapsto \Phi(c)\,,
\end{array}\right.
\end{equation}
we arrive at a stratification $\mathcal{S}=\{S_{j}\}_{j\in J}$ of
$D$ such that for every index $j\in J$,
\begin{itemize}
\item
$\Phi_{j}:=\Phi|_{S_{j}}$ is a $C^{2}$ function of constant rank;
\item
$\Phi_{j}(S_{j})$ is a manifold of dimension equal to the rank of $\Phi_{j}$;\smallskip
\item
the image strata $\{\Phi(S_{j})\}_{j}$ belong to a Whitney
stratification of $\mathbb{R}^{n}$.
\end{itemize}
In particular,
\begin{equation}
D=\bigcup_{j\in J}S_{j}\label{c1}
\end{equation}
and
\begin{equation}
j_{1}\neq j_{2}\quad\Rightarrow\qquad\Phi(S_{j_{1}})=\Phi(S_{j_{2}}
)\quad\text{or}\quad\Phi(S_{j_{1}})\cap\Phi(S_{j_{2}})=\emptyset.\label{c2}
\end{equation}
Denote the set of strata  of full dimension by
$\{S_{j_1},...,S_{j_l}\}$. Set
\[
U=\bigcup_{i=1}^{\ell}S_{j_i}
\]
and observe that the above set is open and dense in $D$, and hence
generic in $S^{n-1}$. \smallskip

Our immediate objective is to show that for every vector $c\in U$
there exists a manifold $\mathcal{M}\subset F$ containing $\Phi(c)$
such that $F$ is partly smooth at $\Phi(c)$ with respect to
$\mathcal{M}$. \smallskip

To this end, fix $\bar{x}\in\Phi(U)$ and consider the set of
``active'' indices
\begin{equation}
I(x) ~:=~ \{j\in J:\bar{x}\in\Phi(S_{j})\}\,.\label{a3}
\end{equation}
We aim to show that the set $F$ is partly smooth at $\bar{x}$
relative to the manifold
\begin{equation}
\mathcal{M}=\Phi_{j}(S_{j}), \label{a4}
\end{equation}
for any $j \in I(x)$. Note that in view of property (\ref{c2}) the
definition of $\mathcal{M}$ is in fact independent of the choice of
$j$ in $I(x)$, and for the same reason the set of active indices
$I(x)$ is invariant for all $x\in\mathcal{M}$. In the sequel, this
set will be simply denoted by $I$.
\smallskip

Clearly, property (i) of the definition of partial smoothness
(Definition~\ref{Definition_partial-smoothness}) holds. If we can
prove properties (ii) and (iii), then our result will follow: since
$U\subset D$, Corollary~\ref{dense} (Density of functionals with
strong maximizers) implies strong criticality for any objective $c
\in U$.

\subsubsection*{Step 1:  normal cone continuity}
We establish the
continuity at $\bar{x}$ of the normal cone mapping $x \mapsto N_{F}(x)$ as $x$
moves along the manifold $\mathcal{M}$.

The normal cone mapping is always outer semicontinuous (even in
$F$). To establish that the truncated normal cone mapping
\begin{equation}
x \mapsto \tilde{\Phi}^{-1}(x)= N_{F}(x)\cap S^{n-1}\qquad (x\in
\mathcal{M})
\label{a5}
\end{equation}
is inner semicontinuous (which clearly suffices for our purposes),
we decompose the above mapping with respect to the active strata. We
set
\begin{equation}
N_{j}(x)=N_{F}(x)\cap S_{j},\qquad\text{for every }j\in J.\label{a6}
\end{equation}
Note that for each $x\in\mathcal{M}$ we have
\begin{equation}
N_{j}(x) \neq \emptyset
~\Leftrightarrow~
j\in I
~\Leftrightarrow~
\mathcal{M}=\Phi(S_{j}).\label{a7}
\end{equation}
We can therefore decompose the truncated normal cone mapping (\ref{a5}) as follows:
\begin{equation}
N_{F}(x)\cap S^{n-1} ~=~ N_{\ast}(x) \cup \bigcup_{j\in I} N_{j}(x)  \label{a8}
\end{equation}
where
\[
N_{\ast}(x)=N_{F}(\bar{x})\cap N_{\ast},
\]
and the set $N_{\ast}$ is defined by equation (\ref{a2}).

\bigskip
\noindent
{\bf Claim A}. For every $x\in\mathcal{M}$ the set $\cup_{j\in
I}N_{j}(x)$ is dense in $N_{F}(x)\cap S^{n-1}$.

\bigskip
\noindent
{\bf Proof of Claim A}. Since we are assuming $\bar{x}\in\Phi(U)$, there exists
an active index $j_{p}$ with $p\in\{1,\ldots,\ell\}$
corresponding to a full-dimensional stratum $S_{j_{p}}$ such that
$\mathcal{M}=\Phi_{ j_p  }(S_{ j_p })$ (see property (\ref{a7})).
This yields that for every $x\in\mathcal{M}$ there exists $c\in
S_{ j_p }$ with $x=\Phi(c)$. Hence
\[
c ~\in~ N_{F}(x)\cap S_{ j_p } ~=~ N_{ j_p }(x) ~\subset~
\bigcup_{j\in I} N_{j}(x).
\]

Fix now any vector $c_{\ast}\in N_{F}(x)\cap S^{n-1}$, and consider the
spherical path
\[
c_{t}:=\frac{c+t(c_{\ast}-c)}{||c+t(c_{\ast}-c)||},\qquad\text{for }
t\in\lbrack0,1].
\]
It follows that $c_{t}\in\mbox{\rm ri}\,N_{F}(x)$, for all
$t\in\lbrack0,1)$. Since $c\in S_{ j_p }\subset D$, there exists a
constant $\delta_{c}>0$ such that $\langle c,x\rangle\geq\langle
c,x^{\prime}\rangle+\delta _{c}||x-x^{\prime}||^{2}$, for all
$x^{\prime}\in F$. By the definition of the normal cone, we also
have $\langle c_{\ast},x\rangle\geq\langle
c_{\ast},x^{\prime}\rangle$ for all $x^{\prime}\in F$. Multiplying
the aforementioned inequalities by $(1-t)$ and $t$ respectively, and
adding, we infer that $x$ is a strong maximizer of $\langle c_{t},
\cdot \rangle$ over the set $F$ for all $0\leq t<1$. In other words,
$c_{t}\in N_{F}(x)\cap D$, which in view of equation (\ref{c1})
yields $c_{t}\in\cup_{j\in I}N_{j}(x),$ for $t\in\lbrack0,1).$ Since
$c_t \to c_{\ast}$ as $t \uparrow 1$, Claim A follows.
\bigskip

In view of Claim A, it is sufficient to establish the inner continuity of
the mapping
\begin{equation}
x \mapsto \bigcup_{j\in I}N_{j}(x)\qquad
x\in\mathcal{M} \text{.}\label{a9}
\end{equation}
To see this, we use the following simple exercise.

\begin{lemma}\label{Lemma}
Let $X$ and $Y$ be metric spaces, and consider two set-valued mappings
$G,T:X\rightrightarrows Y$ such that
$\mbox{\rm cl}(G(x)) = T(x)$ for all points $x\in X.$
If $G$ is inner semicontinuous at a point $\bar{x} \in X$, then so is $T$.
\end{lemma}

\noindent {\bf Proof of Lemma~\ref{Lemma}}. Assume (towards a
contradiction) that there exists a constant $\rho>0$, a sequence
$\{x^k\}\subset X$ with $x^k \rightarrow\bar{x}$ and a point
$\bar{y}\in T(\bar{x})$, such that
\[
\mathrm{dist}(\bar{y},T(x^k))>\rho>0.
\]
Then pick any point $\hat{y}\in B(\bar{y},\rho/2)\cap G(\bar{x})$
and use the inner semicontinuity of $G$ to get a sequence $y^{k}\in
G(x^k)\subset T(x^k)$ for $k\in\mathbb{N}$ such that
$y^{k}\rightarrow\hat{y}$. This gives a contradiction, proving the
lemma.  $\hfill\Box$
\bigskip

\noindent
Applying this lemma to the set-valued mappings
\[
G(x) = \bigcup_{j\in I}N_{j}(x)\qquad\text{and}\qquad
T(x)=N_{F} (x)\cap S^{n-1}
\]
accomplishes the reduction we seek.
\bigskip

Let us, therefore, prove the inner semicontinuity of the mapping
defined in (\ref{a9}) at the point $\bar x$. To this end, fix any
vector $\bar{c}\in\cup_{j\in I}N_{j}(\bar x)$ and consider any
sequence $\{x_{k}\}_{k}\subset\mathcal{M}$ approaching $\bar{x}$.
For some index $j\in I$ we have $\bar{c}\in S_{j}$. Let us restrict
our attention to the constant-rank surjective mapping $\Phi_{j}
:S_{j}\rightarrow\mathcal{M}$ and let us recall that
\[
\Phi_{j}(S_{j})=\mathcal{M}\qquad\text{and}\qquad\Phi_{j}(\bar{c})=\bar
{x}.
\]
Let $d$ be the dimension of the stratum $S_{j}$, so
\[
\mathrm{rank\,}(d\Phi_{j})=\dim\mathcal{M}:=r\leq d\leq n-1.
\]
Denote by $0_{d}$ (respectively $0_{r}$) the zero vector of the
space $\mathbb{R}^{d}$ (respectively $\mathbb{R}^{r}$). Then
applying the Rank Theorem (Proposition~\ref{Proposition_constant}),
we infer that for some constants $\delta,\varepsilon>0$ there exist
diffeomorphisms
\begin{equation}
\psi_{1}:B(0_{d},\delta)\rightarrow S_{j_{0}}\cap
B(\bar{c},\varepsilon
)\quad\text{and}\quad\psi_{2}:B(0_{r},\delta)\rightarrow\mathcal{M}\cap
B(\bar{x},\varepsilon) \label{b3}
\end{equation}
such that
\begin{equation}
\psi_{1}(0_{d})=\bar{c}\quad\text{and}\quad\psi_{2}(0_{r})=\bar{x},
\label{b4}
\end{equation}
and such that all vectors $y\in B(0_{d},\delta)$ satisfy
\begin{equation}
(\psi_{2}^{-1}\circ\Phi_{j}\circ\psi_{1})(y)=\pi(y), \label{b5}
\end{equation}
where for $y=(y_{1},\ldots,y_{d})\in\mathbb{R}^{d}$ we have
\begin{equation}
\pi(y_{1},\ldots,y_{r},y_{r+1}\ldots,y_{d}) ~=~ (y_{1},\ldots,y_{r})
~\in~ B(0_{r},\delta) ~\subset~ \mathbb{R}^{r}. \label{b6}
\end{equation}
We may assume $\{x_{k}\}_{k}\subset\mathcal{M}\cap
B(\bar{x},\varepsilon)$. Thus, in view of definition (\ref{b3}), for
every integer $k\in \mathbb{N}$ there exists a vector $z^{k}
=(z_{1}^{k},...,z_{r}^{k}) \in B(0_{r},\delta)$ with
\begin{equation}
\psi_{2}(z^{k})=x_{k}.\label{b}
\end{equation}
Note $z^{k}\rightarrow  0_{r} = (\psi_{2})^{-1}(\bar{x})$.  Define
vectors
\[
y^{k}:=(z_{1}^{k},...,z_{r}^{k},0,..,0)\in\mathbb{R}^{d}
\]
for every $k\in\mathbb{N}$. Since $z^{k}\in B(0_{r},\delta)$, we
know $y^{k}\in B(0_{d},\delta)$, and clearly
\begin{equation}
y^{k} \rightarrow 0_{d}.\label{b7}
\end{equation}
We now define vectors $c_{k}:=\psi_{1}(y^{k})$ for each $k$. In view
of definition (\ref{b3}) we see that $c_{k}\in S_{j}\cap
B(\bar{c},\varepsilon)$, and in view of properties (\ref{b7}) and
(\ref{b4}),
\[
c_{k} \rightarrow \psi_{1}(0_{d}) = \bar{c} ~~\mbox{as}~ k \to
\infty.
\]

To complete the proof of inner semicontinuity, it remains to show
$c_{k}\in N_{F}(x_{k})$. Since $\Phi_{j}
(c_{k})=\Phi_{j}(\psi_{1}(y^{k}))$ we infer by properties (\ref{b5})
and (\ref{b7}) that
\[
\psi_{2}^{-1}(\Phi_{j}(c_{k}))=(\psi_{2}^{-1}\circ\Phi_{j}\circ\psi_{1}
)(y^{k})=\pi(y^{k})=z^{k}.
\]
Using now the fact that $\psi_{2}$ is a diffeomorphism we deduce
from equation (\ref{b}) that
$\Phi_{j}(c_{k})=\psi_{2}(z^{k})=x_{k}$. Thus
$c_{k}\in\Phi_{j}^{-1}(x_{k})\subset N_{F}(x_{k})$ which completes
the proof of inner semicontinuity and hence of Step 1.

\subsubsection*{Step 2:  sharpness}
It remains to verify that condition (iii) of
Definition~\ref{Definition_partial-smoothness}, namely
\begin{equation} \label{conclusion}
N_{\mathcal{M}}(\bar{x})=N_{F}(\bar{x})-N_{F}(\bar{x}),
\end{equation}
is also fulfilled.

To this end, as in the proof of  Claim A, we can choose an
index $ j \in I$ corresponding thus to a stratum $S_{ j }$ of
full dimension (that is, equal to $n-1$) such that
$\mathcal{M}=\Phi_{ j }(S_{ j })$. Recall that the definable $C^{2}$-mapping
$\Phi_{ j }:S_{ j }\rightarrow\mathcal{M}$ is surjective and has
constant rank $r=\dim\mathcal{M}$, so
\[
\dim N_{\mathcal{M}}(\bar{x})=n-r.
\]
It follows directly from the inclusion $\mathcal{M}\subset F$ that
$N_{F}(\bar{x}) \subset N_{\mathcal{M}}(\bar{x})$.  Since the right-hand side is a subspace, we deduce
\begin{equation}
N_{F}(\bar{x}) - N_{F}(\bar{x}) \subset
N_{\mathcal{M}}(\bar{x}).\label{d1}
\end{equation}
Since $\Phi_{ j }$ is surjective and of maximal rank, we deduce easily that
$\Phi_{ j }^{-1}(\bar{x})$ is a definable submanifold of $S^{n-1}$ of
dimension
\[
\dim\,\Phi_{ j }^{-1}(\bar{x})=(n-1)-r.
\]
which, in view of definition (\ref{a5}) and equation (\ref{a6}) yields
\[
\dim\,\left(  N_{F}(\bar{x})\cap S^{n-1}\right)
\geq\dim\,N_{ j }(\bar {x})\geq(n-1)-r\text{.}
\]
Thus $\dim\,N_{F}(\bar{x})\geq n-r,$ which, along with inclusion (\ref{d1}), yields
equation (\ref{conclusion}), as required. $\hfill\Box$

\bigskip

It is interesting to revisit the example in the introduction,
problem (\ref{bad example}). One can truncate the feasible region
(by intersecting with the unit ball for example) to obtain a convex
compact semi-algebraic set for which the functional
$\bar{c}=(0,0,-1)$ has a unique maximizer (the origin) while the
generic condition of Theorem \ref{Theorem_main} fails. In other
words, $\bar{c} \not\in U$ according to notation of Theorem
\ref{Theorem_main}.
\bigskip

\noindent\textbf{Acknowledgment} A major part of this work was
accomplished during a research visit of the second and third authors
at the University Pierre et Marie Curie (Paris~6), in May 2008.

\begin{center}
----------------------------------------------------
\end{center}

\bigskip

\noindent J\'{e}r\^{o}me BOLTE

\smallskip

\noindent UPMC Univ Paris 06, \'{E}quipe Combinatoire et Optimisation (UMR
7090), Case 189\newline Universit\'{e} Pierre et Marie Curie, 4 Place Jussieu,
F--75252 Paris Cedex 05.

\textit{and}\newline INRIA Saclay, CMAP, Ecole Polytechnique, F--91128
Palaiseau, France.

\smallskip

\noindent E-mail: \texttt{bolte@math.jussieu.fr }; \\
\texttt{http://www.ecp6.jussieu.fr/pageperso/bolte}

\smallskip

\noindent Research supported by the ANR grant ANR-05-BLAN-0248-01 (France).

\vspace{0.7cm}

\noindent Aris DANIILIDIS\smallskip

\noindent Departament de Matem\`{a}tiques, C1/320\newline Universitat
Aut\`{o}noma de Barcelona\newline E--08193 Bellaterra (Cerdanyola del
Vall\`{e}s), Spain.

\textit{and}\newline Laboratoire de Math\'{e}matiques et Physique
Th\'{e}orique\newline Universit\'{e} Fran\c{c}ois Rabelais, Facult\'{e} des
Sciences\newline Parc de Grandmont, F--37200 Tours, France.

\smallskip

\noindent E-mail: \thinspace\texttt{arisd@mat.uab.es}\quad;\quad
\texttt{http://mat.uab.es/\symbol{126}arisd}

\smallskip

\noindent Research supported by the grant MTM2008-06695-C03-03/MTM
(Spain).

\vspace{0.7cm}

\noindent Adrian LEWIS\smallskip

\noindent School of Operations Research and Information Engineering\newline
Cornell University\newline Ithaca, NY 14853, USA.\newline\medskip\noindent
E-mail: \texttt{aslewis@orie.cornell.edu} \\
\texttt{http://www.orie.cornell.edu/\symbol{126}aslewis}\smallskip

\noindent Research supported in part by National Science Foundation Grant
DMS-0806057 (USA).


\begin{thebibliography}{99}

\bibitem {BDL04}\textsc{Bolte, J., Daniilidis, A. \& Lewis, A.S.,} The
\L ojasiewicz inequality for nonsmooth subanalytic functions with applications
to subgradient dynamical systems, \textit{SIAM J. Optim.} \textbf{17} (2007), 1205--1223.

\bibitem {BDL2008}\textsc{Bolte, J., Daniilidis, A. \& Lewis, A.S.,} Tame
functions are semismooth, \textit{Math. Program}. (\textit{Series B})
\textbf{117} (2009), 5--19.

\bibitem {Bur88}\textsc{Burke, J.V. \& J.J. Mor\'e,} On the identification of active constraints, \textit{SIAM J. Numerical Analysis}.
\textbf{25} (1988), 1197--1211.

\bibitem {Coste99}\textsc{Coste, M.,} \textit{An Introduction to o-minimal
Geometry}, RAAG Notes, 81 pages, Institut de Recherche Math\'{e}matiques de
Rennes, November 1999.

\bibitem {DonZol93}\textsc{Dontchev, A.L. \& Zolezzi, T.,} \textit{Well-Posed
Optimization Problems,} Lecture Notes in Mathematics \textbf{1543} (Springer
Verlag, Berlin 1993).

\bibitem {Dries-Miller96}\textsc{van den Dries, L. \& Miller, C.,} Geometric
categories and o-minimal structures, \textit{Duke Math. J.} \textbf{84}
(1996), 497--540.

\bibitem {Ioffe2008}\textsc{Ioffe, A.D., }An invitation to Tame Optimization, 23
p., preprint 2008.

\bibitem {IoffeLuc2005}\textsc{Ioffe, A.D. \& Lucchetti, R., }Typical convex
program is very well posed, \textit{Math. Program. }\textbf{104} (2005), 483--499.

\bibitem {Kurdyka98}\textsc{Kurdyka, K., }On gradients of functions definable
in o-minimal structures, \textit{Ann. Inst. Fourier} \textbf{48} (1998), no.
3, 769--783.

\bibitem {Lee2003}\textsc{Lee, J.M.}, \textit{Introduction to Smooth Manifolds},
(Springer, New York 2003).

\bibitem {LOS2000}\textsc{Lemarechal, C., Oustry, F. \& Sagastizabal, C., }The
$\mathcal{U}$--Lagrangian of a convex function, \textit{Trans. Amer. Math.
Soc.} \textbf{352} (2000), 711--729.

\bibitem {Lewis2002}\textsc{Lewis, A.S.,} Active sets, nonsmoothness and
sensitivity, \textit{SIAM Journal on Optimization} \textbf{13} (2003), 702--725.

\bibitem {Lewis2004}\textsc{Hare, W.L. \& Lewis, A.S.,} Identifying active constraints via partial smoothness and prox-regularity, \textit{Journal of Convex Analysis} \textbf{11} (2004), 251--266.

\bibitem {Phelps93}\textsc{Phelps, R.R.}, \textit{Convex Functions, Monotone
Operators and Differentiability}, Lecture Notes in Mathematics \textbf{1364}
(2nd ed.) (Springer Verlag, Berlin 1993).

\bibitem {Rock98}\textsc{Rockafellar, R.T. \& Wets, R.J.-B.,} \textit{Variational
Analysis}, Grund\-lehren der Mathematischen, Wissenschaften, Vol. \textbf{317} ,
(Springer, 1998).

\bibitem {TaLeLoi96}\textsc{Loi, Ta L\^{e},} Thom stratifications for
functions definable in o-minimal structures on ($\mathbb{R}$,+,$\cdot$),
\textit{C. R. Acad. Sci. Paris S\'{e}r. I Math}. \textbf{324} (1997), no. 12, 1391--1394.

\bibitem {Wri93}\textsc{Wright, S.J.,} Identifiable surfaces in constrained optimization, \textit{SIAM Journal on Control and Optimization} \textbf{31} (1993), 1063--1079.

\end{thebibliography}
\end{document}